\documentclass[12pt]{article}
\usepackage{amsmath,amssymb,amsfonts,amsthm,amscd}
\def\Q{\Bbb Q}

\def\F5{\Bbb F_5}

\def\Z{\Bbb Z}

\def\C{\Bbb C}

\def\R{\Bbb R}
\def\e2{{\eta_2}}

\def\no{\noindent}

\def\s{\subset}

\def\ti{\times}

\def\k{$k$}

\def\op{\oplus}

\def\part{\partial}

\def\sk{{\rm Spec}k}
\def\con{{\rm cone}}
\def\oY{\overline{Y}}
\def\oX{\overline{X}}
\def\oU{\overline{U}}
\def\oH{\overline{H}}
\begin{document}

\newtheorem{theorem}{Theorem}[section]
\newtheorem{proposition}{Proposition}[section]
\newtheorem{lemma}{Lemma}[section]
\newtheorem{corollary}{Corollary}[section]

\title{Nori's construction and the second Abel-Jacobi map}
\author{Kenichiro Kimura}
\date{}

\maketitle 
\section{Introduction}

Let $k$ be a subfield of $\C$. Nori constructs an abelian category of
mixed motives over $k$. One of the fundamental facts in his construction
is the following(\cite{No}):

\begin{theorem}[Basic Lemma]
Let $X$ be an affine scheme of finite type over $k$.  Let $n$ be the
dimension of $X$.  Let $F$ be a weakly constructible
sheaf on $X(\C)$ for the usual topology. 
Then there is a Zariski open $U$ in $X$ with the properties below,
where $j:\,\,U\to X$ denotes the inclusion.  
\begin{enumerate}

 \item dim$Y< n$ where $Y=X-U$.

 \item $H^q(X(\C),j_!j^*F)=0$ for $q\neq n$.

\end{enumerate}
\end{theorem}

Here a sheaf $F$ on $X$ is weakly constructible if $X$ is the disjoint union of finite 
collection of locally closed subschemes $Y_i$ defined over $k$ such that the restrictions
$F|_{Y_i}$ are locally constant.
Beilinson (\cite{Be}, Lemma 3.3) proves this fact in all characteristics
of the base field.
Based on Theorem 1.1 Nori shows that affine $k$-varieties have 
a kind of ``cellular decomposition''. In section 2
of this note we  give an exposition of the outcome if we apply
this construction to etale cohomology. It can also be viewed as a partial
exposition of $l$-adic realization of Nori's category of motives.
The main result (Theorem 2.2)
says that $Rf_*\Z_l(a)$ for a variety
$X\overset{f}\to \sk$ is a complex each component of which comes from
a mixed motive. This gives an answer to a question asked by Jannsen
in a remark in \cite{Ja2}. In \cite{Be} a similar result
for perverse sheaves is proven. 
We learned about Nori's category in \cite{Le}.

\no In section 3, we give a simple description of the second $l$-adic
Abel-Jacobi
map for certain algebraic cycles on a smooth projective variety. We briefly
recall the definition of the $l$-adic Abel-Jacobi map.
Let $X\overset{f}\to \sk$ be a smooth projective variety of dimension $n$.
We denote the absolute Galois
group of $k$ by $G_k$. For an algebraic cycle $z$ on $X$
of codimension $i$ the class

\[[z]\in H^{2i}_{\rm cont}(X,\Z_l(i))\] is defined. Here
$H^{2i}_{\rm cont}(X,\Z_l(i))$ is continuous etale cohomology.
The usual cycle class $cl(z)=cl^0(z)$ is the image of
$[z]$ under the edge homomorphism
\[H^{2i}_{\rm cont}(X,\Z_l(i))\to H^{2i}(\oX,\Z_l(i)).\]

\no The Hochschield-Serre spectral sequence

\[E_2^{p,q}=H^p(G_k,H^q(\oX,\Z_l(i)))\Longrightarrow H^{2i}_{\rm cont}(X,\Z_l(i))\]
induces higher classes
\[cl^j:\,{\rm Ker}cl^{j-1}\to H^j(G_k,H^{2i-j}(\oX,\Q_l(i))).\]
We refer the reader to \cite{Ja2} for more details. See also
\cite{Ra}.

\no Let $z\in
CH^i(X)$ be an algebraic cycle such that\\
$0=cl^1(z)\in H^1(G_k,H^{2i-1}(\oX,\Z_l(i))).$
In Theorem 3.1  we give a simple
description of the push-out of $cl^2(z)\in H^2(G_k,H^{2i-2}(\oX,{\Q_l}(i)))$
by a quotient map $H^{2i-2}(\oX,{\Q_l}(i))\to H^{2i-2}(\oX,{\Q_l}(i))/
 H^{2i-2}_H(\oX, \Q_l(i))$ for a certain multiple hypersurface
section $H$ of $X$.  

\no G. Welters(\cite{We}) also gives a description of the second Abel-Jacobi
map for zero-cycles. It would be interesting to compare these two
descriptions.

\no {\it Acknowledgements}

\no Theorem 3.1 used to be only about zero cycles. 
The author is grateful to Spencer Bloch for pointing out the
possibility of generalizing it to cycles of other dimensions.
He also Thanks Masanori Asakura for his comments on an earlier version. 

\section{Nori's construction}

Let $k$ be a subfield of $\C$. In this note a variety is an 
integral separated scheme of 
finite type over \k. All schemes and morphisms between them
are defined over $k$. For a variety $X$ $\overline{X}$ denotes
$X\ti_{{\rm Spec}\,k}{\rm Spec}\,\bar{k}$ where $\bar{k}$ is an algebraic
closure of \k.


\begin{theorem}
Let $X$ be an affine scheme of finite type of 
dimension $n.$ Let $F$ be a constructible
sheaf on $X(\C)$ for the usual topology. 
Then there is a Zariski open $U$ in $X$ with the properties below,
where $j:\,\,U\to X$ denotes the inclusion.  
\begin{enumerate}

 \item dim$Y< n$ where $Y=X-U$.

 \item $H^q(X(\C),j_!j^*F)=0$ for $q\neq n$.

\end{enumerate}
\end{theorem}

\no Let $X'$ be the largest open subset of $X$ such that $X'$ is smooth
and $F|_{X'}$ is locally constant.  As in Remark 1.1 in
\cite{No} the open set $U$ depends only on the open $X'$ and not on $F$.

\no Fix an integer $a$ and a prime $l$. We are going to Use Theorem 2.1
in the case where $F$ is an etale sheaf of the form ${j_V}_!j_V^*\Z/l^m\Z(a)$ for an open set $V\overset{j_V}{\hookrightarrow}
X$. In this case $H^q(X(\C),j_!j^*F)$ in the assertion 2 of
Theorem 2.1 is isomorphic
to $H^q(\oX_{et}, j_!j^*F)$. 


\no Let $X$ be an affine variety and let $Z$ be a proper closed subset
of $X$. Let $j$ be the inclusion $X-Z\to X.$ 
By applying Theorem 2.1 to $F= 
j_!j^*\Z/l^m\Z(a)$ we obtain the following:
\begin{corollary}
Let $X$ be an affine variety of dimension $n$ and let $Z$ be a closed subset
of dimension$<n$.
Then there exists a closed subset $Y$ of $X$ of dimension$<n$ which contains $Z$ such that

\[H^q(\oX,\oY,\Z/l^m\Z(a))=0\,\,{\rm for} \quad q \neq n \] for all $m\geq 1.$

\end{corollary}
By Corollary 2.1 there is a filtration by closed subsets

\[\emptyset =X_{-1}\subset X_0\subset \cdots \subset X_n=X\]
such that $H^q(\oX_i,\oX_{i-1},\Z/l^m\Z(a))=0$ for $q\neq i$ and for $m\geq 1.$

\no Let $f:\,X\to {\rm Spec}\,k$ be the structure morphism. We consider
$Rf_*\Z_l(a)$ in $D^b(Sh({\rm Spec}\,k_{et})^{\Z_l})$. Here for a variety
$X$ $Sh(X_{et})^{\Z_l}$ denotes the category of $l$-adic sheaves of
Jannsen(\cite{Ja1}, (6.9)). $Sh(X_{et})^{\Z_l}$ has enough injectives.
\begin{proposition}
Let $X$ be an affine variety of dimension $n$ and take a filtration  
\[\emptyset =X_{-1}\subset X_0\subset \cdots \subset X_n=X\] by closed 
subsets as above.

\begin{enumerate}

\item Let $D_n$ be the complex

\[\begin{split}
0&\to H^0(\oX_0,\Z_l(a))\to\cdots \to H^i(\oX_i,\oX_{i-1},\Z_l(a))
\to \cdots \to \\
&H^n(\oX,\oX_{n-1},\Z_l(a))\to 0
\end{split}\]  where the maps between components are the boundary map of cohomology
multiplied by $-1$.
Let $f:\,X\to {\rm Spec}k$ be the structure morphism. Then 
there is a natural isomorphism

\[\R f_*\Z_l(a)  \simeq D_n\]
in $D^b(Sh({\rm Spec}k)^{\Z_l})$. 

\item Let $Y\overset{j_Y}\hookrightarrow X$ be an affine open 
subset of $X$ and \[\emptyset =Y_{-1}\subset Y_0\subset \cdots \subset Y_n=Y\]
be a filtration by closed subsets 
such that $H^q(\oY_i,\oY_{i-1},\Z/l^m\Z(a))=0$ for $q\neq i$ and for $m\geq 1.$ Assume that
$Y_i\subset X_i$ for each $i$. Let 
${D_n}_Y$ be the complex 
\[\begin{split}
0\to H^0(\oY_0,\Z_l(a))\to&\cdots \to H^i(\oY_i,\oY_{i-1},\Z_l(a))
\to \cdots \\ \to
& H^n(\oY,\oY_{n-1},\Z_l(a))\to 0
\end{split}\]
Let $g=f\circ j_Y:\,Y\to \sk.$
Then the isomorphisms $\R f_*\Z_l(a) \simeq D_n$ and 
$Rg_*\Z_l(a)\simeq {D_n}_Y$ are compatible with the pull-back $j_Y^*.$

\end{enumerate}
\end{proposition}

\no {\it Proof.} 1.   
For $0\leq a\leq n$ let $X_a^o=X_a-X_{a-1}$, $j_a:X_a^o\hookrightarrow X_a$
 and $i_a:X_a\hookrightarrow X.$ Let
$B_a^m=Rf_*{i_a}_*\Z/l^m\Z(a)|_{X_a}$. We have exact triangles
\[Rf_*{i_a}_*({j_a}_!\Z/l^m\Z(a)|_{X_a^o})\to B_a^m\to B_{a-1}^m\to.\]
We denote the 
term on the left end by $A_a^m$. 

\no Let $D_n^m$ be the complex

\[\begin{split}
0\to H^0(\oX_0,\Z/l^m\Z(a))\to\cdots &\to H^i(\oX_i,\oX_{i-1},\Z/l^m\Z(a))
\to \\ &\cdots \to 
H^n(\oX,\oX_{n-1},\Z/l^m\Z(a))\to 0
\end{split}\]
\no We need to construct a chain of quasi-isomorphisms between $D_n^m$ and $B_n^m$.
We fix notation and the signs. The mapping cone of $f:A^\cdot \to B^\cdot$
between two cohomological complexes is given in degree $i$ by

\[A^{i+1}\op B^i,\quad{\rm with \,\,differential}\quad
d(a,b)=(-da,db-f(a)).\]
We have a quasi-isomorphism
\[\con(B^m_{n-1}\to \con(B^m_n\to B^m_{n-1}))[-1]\to B^m_n.\]

\no By smooth base change theorem 
$R^qf_*{j_n}_!\Z/l^m\Z(a)=0$ for $q\neq n$ and\\ $R^nf_*{j_n}_!\Z/l^m\Z(a)=H^n(\oX,\oX_{n-1},\Z/l^m\Z(a)).$
Here we regard an etale sheaf on Spec$k$ as a $G_k$ module.
So \[{\rm cone}(B^m_n\to B^m_{n-1})\overset{qis}{\sim} 
 A_n[1]\overset{qis}{\sim}H^n(\oX,\oX_{n-1},\Z/l^m\Z(a))[1-n].\]
So there is a quasi-isomorphism
\[\begin{split}\con(B^m_{n-1}\to &\con(B^m_n\to B^m_{n-1}))[-1]\to\\
&\con(B^m_{n-1}\to 
H^n(\oX,\oX_{n-1},\Z/l^m\Z(a))[1-n])[-1].
\end{split}\]
We have a quasi-isomorphism 
\[\con(B^m_{n-2}\to \con(B^m_{n-1}\to B^m_{n-2}))[-1]\to B^m_{n-1}.\] So there is a quasi-isomorphism
\[
\begin{split} \Bigl[\con\Bigl(\con\bigl(B^m_{n-2}\to &\con(B^m_{n-1}\to B^m_{n-2})\bigr)[-1]\\
&\to H^n(\oX,\oX_{n-1},\Z/l^m\Z(a))[1-n]\Bigr)[-1]\Bigr]\\ & \to \con(B^m_{n-1}\to H^n(\oX,\oX_{n-1},\Z/l^m\Z(a))[1-n])[-1].
\end{split}\]
We can repeat this process until we get the complex $D_n^m$.
All the maps appearing here are compatible with the transition map
$\Z/l^{m+1}\Z(a)\to \Z/l^m\Z(a)$.





2.  All the maps in the construction of the isomorphism
$\R f_*\Z_l(a) \simeq D_n$ are compatible with the pull-back $j_Y^*.$ \qed

In general we use $\check{C}$eck construction. Let $X$ be a variety of dimension
$n$. Let $I$ be a finite set $\{1,\cdots, s\}$ and 
let $U_i\,(i\in I)$ be a covering of $X$ by affine open subsets.
For any finite sets of indices $i_0,\cdots,i_p\in I$ we denote the intersection
$U_{i_0}\cap \cdots  \cap U_{i_p}$ by $U_{i_0\cdots i_p}.$ We denote the open
immersion $U_{i_0\cdots i_p} \hookrightarrow X$ by $j_{i_0\cdots i_p}$.
For $U_{1\cdots s}=\cap_{i\in I}U_i$
choose a filtration by closed subsets $\emptyset=C_{1 \cdots s}^{-1} \s C_{1 \cdots s}^0
\s C_{1 \cdots s}^1\s \cdots C_{1 \cdots s}^{n-1}\s
C_{1 \cdots s}^n=U_{1 \cdots s}$ such that $H^j(\overline{C_{1\cdots s}^i},\overline{C_{1 \cdots s}^
{i-1}},\Z/l^m\Z(a))=0$ for $i\neq j$ and for $m\geq 1$. When the filtrations 
$C_{i_0\cdots i_p}^{j}$ of $U_{i_0\cdots i_p}$ are defined for all subsets 
$\{i_0,\cdots,i_p\}\s I$ with $p+1$ elements and for $0 \leq j\leq n-1$ then the 
filtrations
$ \emptyset =C_{i_0 \cdots i_{p-1}}^{-1}\s C_{i_0 \cdots i_{p-1}}^0
\s C_{i_0 \cdots i_{p-1}}^1\s \cdots C_{i_0 \cdots i_{p-1}}^{n-1}\s
C_{i_0 \cdots i_{p-1}}^n=U_{i_0 \cdots i_{p-1}}$ are chosen so that
$H^j(\overline{C_{i_0 \cdots i_{p-1}}^i},\overline{C_{i_0 \cdots i_{p-1}}^
{i-1}},\Z/l^m\Z(a))=0$ for $i\neq j$ and for $m\geq 1$ and also that 
$C_{i_0 \cdots i_p}^i \s C_{i_0 \cdots i_{p-1}}^i$ for any subset 
$\{i_0, \cdots ,i_p\}\supset \{i_0, \cdots ,i_{p-1}\}$.

\no We have the $\check{C}$eck complex

\[\begin{split}
0\to \underset{i\in I}\op R(f\circ j_i)_*\Z_l(a)&\to \underset{i_0<i_1}\op
R(f\circ j_{i_0i_1})_*\Z_l(a)\to\\
& \cdots\to R(f\circ j_{1 \cdots s})_*
\Z_l(a)\to 0 
\end{split}\] The total complex associated to this double
complex is isomorphic to $Rf_*\Z_l(a)$ in $D^b(Sh(\sk_{et})^{\Z_l}).$ By the assertion
1 of Proposition 2.1
$R(f\circ j_{i_0\cdots i_p})_*\Z_l(a)$ is isomorphic to 
\[\begin{split}
0\to H^0(\overline{C^0_{i_0\cdots i_p}},\Z_l(a))\to &H^1(\overline{C^1_{i_0\cdots i_p}},
\overline{C^0_{i_0\cdots i_p}},\Z_l(a))\to \\&\cdots \to  H^n(\overline{U_{i_0\cdots i_p}},\overline{C^{n-1}_{i_0\cdots i_p}},\Z_l(a))
\to 0.\end{split}\] We denote this complex by $D_{i_0\cdots i_p}$. 
By the assertion 2 of Proposition 2.1 the isomorphism $D_{i_0\cdots i_p}\to
R(f\circ j_{i_0\cdots i_p})_*\Z_l(a)$ is compatible with the restriction 
$R(f\circ j_{i_0\cdots i_p})_*\Z_l(a) \to
R(f\circ j_{i_0\cdots i_{p+1}})_*\Z_l(a)$.  So we obtain the following.
\begin{theorem}
The complex $Rf_*\Z_l(a)$ is isomorphic to the total complex associated to
the double complex

\[0\to \underset{i\in I}\op D_i\to \underset{i_0<i_1}\op
D_{i_0i_1}\to \cdots\to D_{1 \cdots s}\to 0 \] in $D^b(Sh(\sk_{et})^{\Z_l})$.
\end{theorem}

\section{A simple description of the second $l$-adic Abel-Jacobi map}

Let $X$ be a smooth projective variety of dimension $n$. 
Let $z\in CH^i(X)$ be an algebraic cycle which is homologous
to $0$. 
Let $cl^1(z)$ be the
image of $z$ under the 
$l$-adic Abel-Jacobi map
\[cl^1:\,CH^i(X)_{\rm hom}\to H^1(G_k,H^{2i-1}(\oX,\Z_l(i))).\]

\no Assume further that the cycle $z$ satisfies the following
condition:

\begin{quote}
{\it Let $q=2i-1-n.$ Then there exists a smooth multiple hypersurface
section $H$ of $X$ of codimension $q$ which supports the cycle $z$ 
such that $z$ is homologous to $0$ on $H$.}

\end{quote}

\no For example if $i=n$ then by Proposition 4.8 of \cite{Ja}(see also Proposition 5.3
of \cite{No2}) such a $H$ always exists. Let $|z|$ be the support of $z$ and
let $Y=X-|z|.$ Let $U=X-H$ and let $j:U\subset Y$ be the inclusion..

\no Let $g: Y\to \sk$ be the structure morphism. The 2-extension 

\[ 0\to H^{2i-2}(\oY,\Q_l(i))\to \frac{Rg_*\Q_l(i)^{2i-2}}{{\rm Im}
\part^{2i-3}}\overset{\part^{2i-2}}\to {\rm Ker}\part^{2i-1}\to
H^{2i-1}(\oY,\Q_l(i))\to 0\] is denoted $\chi_{2i-2}(Y)$ in \cite{Ja2}.

\no Assume that $0=cl^1(z)\in H^1(G_k,H^{2i-1}(\oX,\Z_l(i))) .$ Then by
Theorem 1 in \cite{Ja2} the class $-cl^2(z)\in {\rm Ext}_{G_k}^2(\Q_l,
H^{2i-2}(\oY,\Q_l(i)))$ is the pull-back of $\chi_{2i-2}(Y)$ by the splitting
$cl^1(z):\Q_l\to H^{2i-1}(\oY,\Q_l(i)).$

\no Let $\cal C$ be the complex
\[
0\to H^{2i-1}(\oX,\Q_l(i))\to H^{2i-1}(\oY,\Q_l(i))\to H^{2i}_{|z|}(\oX,\Q_l(i))
\to 0\] and let ${\cal C}_ H$ be the complex 

\[0\to H^{2(i-q)-1}(\oH,\Q_l(i-q))\to H^{2(i-q)-1}(\overline{H\cap Y},\Q_l(i-q))\to
H^{2(i-q)}_{|z|}(\oH,\Q_l(i-q))\to 0. \] There is the Gysin map ${i_H}_*:
{\cal C}_H\to \cal C.$ From the definition of $q$
the map  ${i_H}_*:H^{2(i-q)-1}(\oH,\Q_l(i-q))\to H^{2i-1}(\oX,\Q_l(i))$
is surjective by hard Lefschetz theorem. So the cycle class 
$cl^1(z)\in {\rm Ext}^1_{G_k}(\Q_l,H^{2i-1}(\oX,\Q_l(i)))$
is the image of 
$cl^1(z)\in {\rm Ext}^1_{G_k}(\Q_l,H^{2(i-q)-1}(\oH,\Q_l(i-q)))$ under
the Gysin map. It is also equal to the push-out by the quotient
$H^{2(i-q)-1}(\oH,\Q_l(i-q))\to \frac{H^{2(i-q)-1}(\oH,\Q_l(i-q))}{\part H^{2i-2}(\oU, \Q_l(i))}$.
So there is a splitting 
\[cl^1(z):\,\Q_l\to 
\frac{H^{2(i-q)-1}(\overline{H\cap Y},\Q_l(1))}{\part H^{2i-2}(\oU, \Q_l(n))}.\]

\begin{theorem} The push-out of $-cl^2(z)\in
 {\rm Ext}^2_{G_k}(\Q_l,H^{2i-2}(\oY,\Q_l(i)))$ by the quotient 
 $H^{2i-2}(\oY,\Q_l(i))\to \frac{H^{2i-2}(\oY,\Q_l(i))}{H^{2i-2}_{Y-U}(\oY,\Q_l(i))}
 $ is given by the pull-back of the 2-extension
 \[\begin{split}
0&\to \frac{H^{2i-2}(\oY,\Q_l(n))}{H^{2i-2}_{Y-U}(\oY,\Q_l(i))}\to
 H^{2i-2}(\oU,\Q_l(i))\to H^{2i-1}_{Y-U}(\oY,\Q_l(i))\to\\ 
&\frac{H^{2i-1}_{Y-U}(\oY,\Q_l(i))}{\part H^{2i-2}(\oU,\Q_l(i))}\to 0
\end{split}\]
by $cl^1(z):\Q_l\to \frac{H^{2i-1}_{Y-U}(\oY,\Q_l(i))}{\part H^{2i-2}(\oU,\Q_l(i))}.$
\end{theorem}

\no {\it Remark}. 

\no When $i$ is equal to $n$, $H^{2n-2}_{Y-U}(\oY,\Q_l(n))$
is generated by the cohomology class of $H(1)$. So we do not lose too
much information by the push-out.

\quad

{\it Proof.} We have an exact triangle
\[ Rg_*R{i_H}^!\Q_l(i) \overset{{i_H}_*}\to Rg_*\Q_l(i)\overset{j^*}\to
R(g\circ j)_* \Q_l(i)|_U\to.\] We denote this triangle by
$A\overset{{i_H}_*}\to B\overset{j^*}\to C\to .$ 
The 2-extension $\chi_{2i-2}(Y)$ is given by

\[0\to H^{2i-2}(\oY)\to \frac{B^{2i-2}}{\part_B^{2i-3}(B^{2i-3})}
 \overset{\part_B^{2i-2}}\to {\rm Ker}\part_B^{2i-1}\to
H^{2i-1}(\oY)\to 0.\] Let $C_2$ be the complex
\[\begin{split}
0&\to H^{2i-2}(\oY) \to \frac{(j^*)^{-1}{\rm Ker}\part_C^{2i-2}+{i_H}_*(A^{2i-2})}
{\part_B^{2i-3}(B^{2i-3})} \overset{\part^{2i-2}_B}\to\\
&{i_H}_*({\rm Ker}\part^{2i-1}_A)\to \frac{H^{2i-1}_{Y-U}(\oY)}{\part
(H^{2i-2}(\oU))}\to 0
\end{split}\]

\no Let $C_3$ be the complex
\[\begin{split}
0&\to  \frac{H^{2i-2}(\oY)}{{\rm Im}({i_H}_*({\rm Ker}\part_A^{2i-2}))}
\to \frac{(j^*)^{-1}{\rm Ker}\part_C^{2i-2}+{i_H}_*(A^{2i-2})}
{\part_B^{2i-3}(B^{2i-3})+{i_H}_*(A^{2i-2})}\overset{\part^{2i-2}_B}\to \\
&\frac{{i_H}_*({\rm Ker}\part^{2i-1}_A)}{{i_H}_*(\part^{2i-2}_A(A^{2i-2}))}
\to \frac{H^{2i-1}_{Y-U}(\oY)}{\part(H^{2i-2}(\oU))}\to 0
\end{split}\]

\no Let $C_4$ be the complex
\[ 0\to \frac{H^{2i-2}(\oY)}{H^{2i-2}_{Y-U}(\oY)}\to H^{2i-2}(\oU)
\to H^{2i-1}_{Y-U}(\oY) \to \frac{H^{2i-1}_{Y-U}(\oY)}{\part(H^{2i-2}(\oU))}\to 0
\]

\no There are natural maps of complexes

\[\chi_{2i-2}(Y)\leftarrow C_2\to C_3\to C_4.\]  So there are natural maps
between the pull-backs of these complexes 
by the splittings given by $cl^1(z)$. Since $C_2$ is exact
this completes the proof. \qed

Kenichiro Kimura

Institute of Mathematics

University of Tsukuba

Tsukuba, Ibaraki

305-8571

Japan

email: kimurak@math.tsukuba.ac.jp


\begin{thebibliography}{99}
%
\bibitem{Be}Beilinson, A.A.,
On the derived category of perverse sheaves,
{\it $K$-theory, arithmetic and geometry},Lect. Not. in Math. {\bf 1289},
Springer-Verlag,1987, 27--41.
%
\bibitem{Ja1}Jannsen, U.,
{Continuous \'etale cohomology}, Math.Ann. {\bf 280}\,(1988), 207--245.
%
\bibitem{Ja}Jannsen, U.,
Equivalence relations on algebraic cycles,
{\it The Arithmetic and Geometry of Algebraic Cycles (Banff, Alberta, 1998)},
Nato Sci.Ser.C Math.Phys.Sci. {\bf 548}, Kluwer,Dortrecht, 2000, 225--260.


\bibitem{Ja2}Jannsen, U.,
Letter from U. Jannsen to B. Gross on higher Abel-Jacobi maps,
{\it The Arithmetic and Geometry of Algebraic Cycles (Banff, Alberta, 1998)},
Nato Sci.Ser.C Math.Phys.Sci. {\bf 548}, Kluwer,Dortrecht, 2000, 261--275.
%

\bibitem{Le}Levine,M.,
Mixed motives.  {\it Handbook of $K$-theory},Vol.{\bf 1} Springer, Berlin (2005),
 429--521.

\bibitem{No}Nori, M.,
Constructible sheaves, {\it Algebra, arithmetic and geometry,
Part I, II (Mumbai, 2000),}Tata Inst. Fund.Res. Stud. math.,
{\bf 16}, Tata Inst. Fund.Res., Bombay, 2002. 471-491.
%
\bibitem{No2}Nori, M.,
Algebraic cycles and Hodge theoretic connectivity,
Invent. Math. {\bf 111}(1993),349-373.
%
\bibitem{Ra}Raskind, W.,
Higher $l$-adic Abel-Jacobi mappings and filtrations on Chow groups,
Duke Math. J. {\bf 78}(1995), 33-57.

\bibitem{We}Welters, G.,
{The Brauer group and the second Abel-Jacobi map for
0-cycles on algebraic varieties},
Duke Math. J. {\bf 117}\,(2003),447--487.


\end{thebibliography}
\end{document}